\documentclass{article}
\usepackage{latexsym}
\usepackage{amssymb}
\title{\textsf{Centralisers of Involutions in Black Box Groups}}
\author{\textsf{Alexandre V. Borovik}\thanks{The author was supported by
the Royal Society and The Leverhulme Trust.}}
\date{\textsf{21 October 2001}}

\newtheorem{theorem}{Theorem}[section]
\newtheorem{fact}[theorem]{Fact}
\newtheorem{conjecture}[theorem]{Conjecture}
\newtheorem{problem}[theorem]{Problem}
\newtheorem{question}[theorem]{Question}

\begin{document}

\pagestyle{myheadings}
\markright{\scriptsize\textsc{A. V. Borovik $\bullet$
Centralisers of Involutions in Black Box Groups $\bullet$
21.10.01}}

\maketitle

\begin{abstract}
We discuss basic structural properties of finite black box groups.
A special emphasis is made on the use of centralisers of
involutions in probabilistic recognition of black box groups. In
particular, we suggest an algorithm for finding the $p$-core of  a
black box group of odd characteristic. This special role of
involutions suggest that the theory of black box groups
reproduces, at a non-deterministic level, some important features
of the classification of finite simple groups.

2000 Mathematical Subject Classification: 20P05.
\end{abstract}
\section{What is a black box group?}

A {\it black box group} $X$ is a device or an algorithm ({\em
`oracle'} or {\em `black box'\/}) which produces (nearly)
uniformly distributed independent random elements from some
finite group $X$. These elements are encoded as 0--1 strings of
uniform length; given strings representing $x,y \in X$, the black
box can compute strings representing $xy$ and $x^{-1}$, and
decide whether $x=y$ in time bounded from above by a constant. In
this setting, one is usually interested in finding probabilistic
algorithms which allow us to determine, with probability of error
$\epsilon$, the isomorphism type of $X$ in time
$O(|\epsilon|\cdot (\log |X|)^c)$. We say in this
situation that our algorithm is run in \textit{Monte Carlo
polynomial time}. A critical discussion of this concept can be
found in \cite{babai2}, while \cite{babai-beals} contains a
detailed survey of the subject. See also the forthcoming book by
Seress \cite{seress}.

In this paper we discuss a (still rather rudimentary) structural
approach to the theory of black box group. We briefly survey
methods for constructing black box oracles for subgroups and
factor groups of black box groups, and then show how one can
construct black box oracles for centralisers of involutions. They
are used in the algorithm for finding the $p$-core of a
 black box group of characteristic $p$.

Isomorphisms and homomorphisms of black box groups are understood
as isomorphisms and homomorphisms of their underlying groups.
However we reserve the term \textit{black box subgroup} for a
subgroup of a black box group endowed with its own black box
oracle.

Despite this rather abstract general setting, practically
important black box groups usually appear as big permutation or
matrix groups. For example, given two square matrices $x$ and $y$
of size, say, $100$ by $100$ over a finite field, it is
unrealistic to list all elements in the group $X$ generated by
$x$ and $y$ and determine the isomorphism class of $X$ by
inspection. But this can often be done, with an arbitrarily small
probability of error, by studying a sample of random products of
the generators $x$ and $y$. The explosive growth of the theory of
black box groups in recent years is reflected in numerous
publications (see, for example, the survey paper \cite{L-G} on the
computational matrix group project) and the fundamental work
\cite{kantor}), and algorithms implemented in the software
packages GAP \cite{G} and MAGMA \cite{BCP}. Our observation
(Section~\ref{sec:oddcharacteristic}) that centralisers of
involutions allow to compute unipotent radicals of black box
groups of odd characteristic might be used in the computational
matrix group project.

This paper is written by a group theorist, not a probabilist. The
author  had the audacity to list some problems of probabilistic
nature directly related to the computational aspects of the black
box group theory.

\subsection{The order oracle}

Almost nothing can be said about a black box group without access
to additional information. In some cases (for example, when our
black box is given as a permutation group of computationally
feasible degree) we have the so-called \textit{order oracle},
that is, we can determine the orders of elements $x \in {X}$. Of
course, in a permutation group the order of an element can be
easily read from its cycle structure. Another situation when we
can determine the order of an element is when we are given a
reasonably small superset $\pi$ of prime integers dividing the
order $|{X}|$ of ${X}$ as well as reasonable bounds for $|{X}|$.
Then we can make the list of all divisors $d$ of $|{X}|$ and try
all of them  by checking whether $x^d =1$; the minimal such $d$
is, of course, the order of $x$. In the case of matrix groups
${X} \leqslant GL_n(\mathbb{F}_q)$ this means that we have to
factorise $|GL_n(\mathbb{F}_q)|$ into primes, which is as hard as
the general factorisation problem \cite{babai-beals}.

However, there is a satisfactory, for the purpose of practical
computation, way around the problem. Instead of the precise
factorisation of $|GL_n(\mathbb{F}_q)|$ one can use a substitute,
the finest factorisation one can get:
\begin{eqnarray*}
|GL_n({\mathbb F}_q)| & = & q^{n(n-1)/2}(q-1)(q^2-1)\cdots (q^n-1)\\[1ex]
& = & r_1^{a_1} \cdots r_m^{a_m}
\end{eqnarray*}
Now we can use the  \textit{pseudo-order\/} $|x|_{\approx}$ of the
element $x$ instead of its exact order:
\begin{eqnarray*}
|x|_{\approx} & = & \hbox{the least number } l = r_1^{b_1}\cdots
r_m^{b_m} \hbox{ such that } x^l =1.
\end{eqnarray*}

Some of the constructions in the present paper can be carried out
under a milder assumption that we know a \textit{global
exponent}, a number $E$ such that $x^E$ can be easily computed
(the standard square-and-multiply method requires about $2 \log_2
E$ multiplications) and $x^E = 1$ for all  $x \in {X}$.

\subsection{Three types of problems}

Given a black box group ${ X}$, we usually deal with one of the
following problems.

\begin{itemize}
\item \textsl{Identification problem.} Determine the isomorphism type of
$\,{ X}$ with the given degree of certainty.

\item \textsl{Verification problem.} Is $ X$ isomorphic to the given
group $G$?

\item \textsl{Constructive recognition.} Find an explicit isomorphism
${ X} \longrightarrow { G}$.
\end{itemize}

 {}From the probabilistic point of view, these three
problems have different nature. In the identification problem, we
have to be prepared that our algorithm might produce a false
answer, although the probability of this outcome can be made
arbitrarily small by running the algorithm sufficiently many
times. The algorithms for the verification problem are usually
one-sided: for example, if we have found in ${ X}$ an element of
order not present in $ G$ then we \textit{definitely} know that
${ X}$ is not isomorphic to $ G$. Constructive recognition
algorithms are also of probabilistic nature and allow for some
probability of failure. But once succeeded, they provide a proof
of the isomorphism ${ X} \simeq { G}$.

\subsection{The use of involutions in recognition of black
box groups: an elementary example}

The classical  Miller-Rabin primality test (from computational
number theory \cite{rabin}; see also \cite[Section~V.1]{koblitz})
is based on the fact that an odd integer $n$ is prime if and only
if $({\mathbb Z}/n{\mathbb Z})^*$ is the cyclic group of order
$n-1$. The group ${ X} = ({\mathbb Z}/n{\mathbb Z})^*$ is a
nicest possible black box: using standard random number
generators, we can produce random  uniformly distributed
independent elements from ${ X}$. Therefore we are in the
setting  of the verification problem: is ${ X}\simeq
{\mathbb{Z}}_{n-1}$? We work on this assumption, in particular,
we assume that ${ X}$ has a global exponent $E=n-1$, and wish to
detect structural differences between  ${ X}$ and
${\mathbb{Z}}_{n-1}$. Notice that the group ${{\mathbb Z}}_{n-1}$
contains exactly one involution (that is, an element of order
$2$). On the other hand, if $n = p_1^{l_1} \cdots p_k^{l_k}$ is
the prime factorisation of $n$ then
$$
({\mathbb Z}/n{\mathbb Z})^*  =  ({\mathbb Z}/p_1^{l_1}{\mathbb
Z})^* \times \cdots \times ({\mathbb Z}/p_k^{l_k}{\mathbb Z})^*
$$
The Chinese Remainder Theorem allows us to lift the involutions
$-1\, {\rm mod}\, p_i^{n_i}$  to involutions in $({\mathbb
Z}/n{\mathbb Z})^*$, thus showing that the involutions in
$({\mathbb Z}/n{\mathbb Z})^*$ generate an elementary abelian
subgroup of order $2^k$.

The key point of the story is that we can easily compute
involutions in ${\mathbb Z}_{n-1}$, and this simple trick will be
used later in this paper. Indeed, we can factorise $n-1$ into a
power of $2$ and an odd factor: $n-1 = 2^l \cdot m$, $m$ odd.
Obviously, at least half of the elements in ${\mathbb Z_{n-1}}$
are of even order, so, with probability at least $1/2$, $x^m$ is
a non-trivial $2$-element. The last non-identity element in the
sequence of squares
$$
x^m, (x^m)^2, \ldots, (x^m)^{2^l} $$ has order $2$; we denote it
${\rm i}(x)$ and call the \textit{involution produced by} $x$.
For the sake of completeness of this definition, we set ${\rm
i}(x) =1$ if $x$ is of odd order (and thus ${\rm i}(x) =x^m$).

If $({\mathbb Z}/n{\mathbb Z})^* \not\simeq {\mathbb Z}_{n-1}$,
this procedure is likely to fail (that is, ${\rm i}(x)$ is not an
involution), due to the fact that, for most integers $n$,
$(x^m)^{2^l} \ne 1 \bmod n$ with probability at least $1/2$. In
the worst case scenario (when $n$ is a so-called
\textit{Carmichael\/} number), the probability of producing $\pm
1$ can be shown to be less than $\frac{1}{2^{k-1}} \leqslant
\frac{1}{4}$. Hence we come to the following formulation of the
Miller-Rabin primality test.

\bigskip

\noindent \textsc{repeat} for random $x \in ({\mathbb
Z}/n{\mathbb Z})^*$:
\begin{itemize}
\item \textsc{compute} ${\rm i}(x)$.
\item \textsc{if} the computation of ${\rm i}(x)$ fails or  ${\rm i}(x)\ne \pm 1$,
\textsc{return}
\begin{quote}
{\sl $n$ is not prime}
\end{quote}
\item \textsc{if}  ${\rm i}(x)= \pm 1$ for $l$ random values of $x$, \textsc{return}
\begin{quote}
{\sl $n$ is  prime with probability of error $\leqslant
\frac{1}{4^l}$.}
\end{quote}
\end{itemize}

We discuss later in the paper the use of involutions and
centralisers of involutions in the analysis of  black box groups,
and the reader will be likely to agree that our approach can be
viewed as a non-commutative version of the Miller-Rabin primality
test. The role of involutions in identification of simple finite
groups, black box or not, is not surprising to a finite group
theorist, since this is the main tool of the classification of
finite simple groups. But even in the very elementary, from the
group-theoretic point of view, setting of the Miller-Rabin test,
involutions are the keys to the structure of the group. Indeed,
if we know an involution $x \ne -1$ in $({\mathbb Z}/n{\mathbb
Z})^*$, then, since $x^2 - 1 \equiv  0 \bmod n$, we have $n \mid
(x-1)(x+1)$, and the calculation of ${\rm gcd}(n, x-1)$ and ${\rm
gcd}(n, x+1)$ yields a non-trivial factor of $n$. Hence the
knowledge of involutions in $({\mathbb Z}/n{\mathbb Z})^*$
amounts to factorisation of $n$ into prime numbers. This simple
observation is of considerable practical value since it makes the
basis of Simmons' attack on the repeated use of the modulus in
the RSA encryption algorithm \cite{boneh}.

\section{Black box groups of odd characteristic}
\label{sec:oddcharacteristic}

A black box group ${X}$ is said to be of (known) characteristic
$p$ if it is isomorphic to a section in the matrix group
$GL_n({\mathbb{F}}_{p^k})$ and the order
$|GL_n({\mathbb{F}}_{p^k})|$ is of computationally feasible size.
Notice that this means, in particular, that we can take $E =
|GL_n({\mathbb{F}}_{p^k})|$ for a global exponent for ${X}$.

A detailed discussion of the following fundamental result can be
found in Babai and Shalev  \cite{bs}. It summarises the work of
\cite{bkps} based on  \cite{kantor} and \cite{altseimer-borovik,borovik}.

\begin{fact}
Given a black box group\/ $X$ of known characteristic, the
standard names of all non-abelian factors of\/ $X$ can be computed
in Monte Carlo polynomial time. \label{fact:names}
\end{fact}

However, the determination of the $p$-core of ${X}$ is an open
problem. Recall that the \textit{$p$-core} $O_p(X)$ is the maximal
normal $p$-subgroup of $X$.

\begin{problem}[{\rm  \cite[Problem~10.2]{babai-beals},
 \cite[Section~5]{bs}}]
\textsf{ Given a black box group ${X}$ of characteristic $p$,
 can one decide  whether $O_p({{X}}) = 1$ in polynomial time?}
 \label{prob:core}
\end{problem}
The answer is not known even if  $O_p({{X}})$  is known to be a
minimal normal subgroup in ${{X}}$ and ${{X}}/O_p({{X}})$ a
simple group of Lie type in characteristic $p$. Moreover, as
shown in \cite{bs}, the general problem can be reduced to this
minimal configuration.

However, in odd characteristic the question can be answered with
the help of centralisers of involutions. Recall that a finite
group $G$ is \textit{quasisimple} if $G = G'$ and $G/Z(G)$ is a
simple group. A \textit{semisimple} group $G$ is a central product
$L_1\cdots L_k$ of  quasisimple groups $L_i$, called {\em
components} of $G$. A \textit{reductive} group $G$ is a central
product of a semisimple group and an abelian $p'$-group. This is,
of course, a finite group theoretic version of the concept of
reductive algebraic group. We say that a reductive group $G$ is
of {\em characteristic $p$}, if all components $L_i$ of $G$ are
quasisimple groups of Lie type.

 For a finite group $H$,
$O^{2'}(H)$ is the subgroup in $H$ generated by all elements of
odd order.

The following theorem puts Problem~\ref{prob:core} (with
$X/O_P(X)$ simple group of Lie type and characteristic $p$) in an
inductive setting.

\begin{theorem}
Let\/ ${X}$ be a black box group of known odd characteristic
$p\geqslant 5$. Assume that\/ $\overline{X}={{X}}/O_p({{X}})$ is a
reductive group of characteristic $p$.  Then we can determine, in
polynomial time, whether\/ $O_p({{X}}) \ne 1$, and, if\/
$O_p({{X}}) \ne 1$, find a non-trivial element from $O_p({{X}})$.
\label{th:Op(X)}
\end{theorem}

An analogous result, although much more technical, can be proven
in characteristic $p=3$; some care is needed in this special case
because of the solvability of small groups like
$SL_2({\mathbb{F}}_3)$.

A similar result was announced by C.~Parker and R.~Wilson.

The proof of Theorem~\ref{th:Op(X)} will be published elsewhere.
Its main idea is to reduce the problem of detecting the
non-trivial $p$-core $O_P(X)$ in $X$ to the similar problem for
the centraliser  of an involution $C_X(t)$. Thus the algorithm is
recursive. Fortunately, the properties of $\overline{X}$ are
inherited by the consecutive centralisers of involutions (or,
what is the same, by the centralisers of elementary abelian
$2$-subgroups) because of the following well-known result on
centralisers of abelian subgroups of semisimple elements
\cite{ste1}.

\begin{fact}
Let\/ $X$ be a reductive group of characteristic $p >5$ and\/
$A<X$ an elementary abelian $2$-group. Set\/ $Y = C_X(A)$. Then
$Y$ contains a reductive normal subgroup $Y^\circ$ of
characteristic $p$ such that $Y/Y^\circ$ is an elementary abelian
$2$-group.
\label{fact:centralisers-of-involutions}
\end{fact}

Therefore we need efficient methods of computation of centralisers
of involutions and various normal subgroups in black box groups.
Because of the probabilistic nature of the algorithm, it is vital
to avoid the possible accumulation of errors.

In this paper, we concentrate on discussion of various problems
related to handling the centralisers of involutions in black box
groups.

\section{Subgroups of black box groups}

A problem which we immediately encounter when dealing with black
box groups is how to construct a good black box for the subgroup
generated by given elements. For example, given a group generated
by a collection of matrices,
$$
X  \leqslant  GL_N({\mathbb F}_q),\quad
 X  =  \langle x_1, \ldots, x_k \rangle,
$$
how can we produce (almost) uniformly distributed independent
random elements from $X$? The commonly used solution is the {\em
product replacement algorithm} \cite{celler2}.

\subsection{The Product Replacement Algorithm}

Denote by $\Gamma_k({X})$ the graph whose vertices are generating
$k$-tuples of elements in ${X}$ and edges are given by the
following transformations:
\begin{eqnarray*}
(x_1, \ldots, x_i, \ldots, x_k) & \longrightarrow & (x_1, \ldots, x_j^{\pm 1}x_i, \ldots, x_k)\\[1ex]
(x_1, \ldots, x_i, \ldots, x_k)  & \longrightarrow & (x_1, \ldots,
x_ix_j^{\pm 1}, \ldots, x_k)
\end{eqnarray*}
The recipe for production of random elements from $ X$ is
deceptively simple: walk randomly over this graph and select
random components $x_i$. The detailed discussion of theoretical
aspects of this algorithm can be found in Igor Pak's survey
\cite{P2}. Pak \cite{P3} has also shown  that, if $k$ is
sufficiently big, the mixing time for a random walk on
$\Gamma_k({X})$ is polynomial in $k$ and $\log |{X}|$. Here, the
{\em mixing time} $t_{\rm mix}$ for a random walk on a graph
$\Gamma$ is the minimal number of steps such that after these
steps
$$
\frac{1}{2} \sum_{v\in \Gamma} \left| P(\hbox{get at } v) -
\frac{1}{\#\Gamma}\right| < \frac{1}{e}.
$$
At the intuitive level, this means that the distribution of
the end points of random walks on $\Gamma$
is sufficiently close to the uniform distribution.

The graph $\Gamma_k({X})$ is still a very mysterious object.
Notice, in particular, that, in general, it is not connected. The
following very natural question is still open.

\begin{conjecture} \textsf{If\/ $G$ is a finite simple group,
the graph\/
$\Gamma_k(G)$ is connected for\/ $k \geqslant 3$.}
\end{conjecture}

However, Pak \cite{P1} found a sufficiently good approximation to
the connectivity of $\Gamma_k(G)$: if $\{\,G_i \mid i =
1,2,\ldots \}$ is a sequence of simple group of increasing order
then one of the connected components of $\Gamma_k(G_i)$ is
asymptotically of the same size as $\Gamma_k(G_i)$.

A remarkable observation by
Lubotzky  and Pak gives a
conceptual explanation of the good properties of the product
replacement algorithm.

\begin{fact}[{\rm Lubotzky  and Pak \cite{LuP}}]
If\/ ${\rm Aut} \, F_k$ satisfies the Kazhdan property {\rm(}T{\rm )},
then mixing time $t_{\rm mix}$ of a random walk on a
component of\/ $\Gamma_k(G)$ is bounded as
$$
t_{\rm mix} \leqslant C(k) \cdot \log_2 |G|.
$$
\end{fact}

Thus the issue is reduced to the long standing conjecture:

\begin{conjecture}
\textsf{For\/ $k \geqslant 4$,  ${\rm Aut} \, F_k$ satisfies the Kazhdan property {\rm
(}T{\rm )}.}
\end{conjecture}

Following \cite{K},
we say that a topological group $G$ satisfies the Kazhdan (T)-property if, for some
compact set  $Q \subset G$,
$$
K = \inf_{\rho} \inf_{v \ne 0} \max_{q \in Q}
\frac{\|\rho(q)(v)-v\|}{\|v\|} > 0,
$$
where $\rho$ runs over all unitary representations of $G$ without
fixed non-zero vectors.
In our context,  ${\rm Aut} \, F_k$ is endowed with the discrete topology.

\section{Normal subgroups}
\label{sec:normal}

Given elements $y_1,\ldots, y_k$ of a black box group ${ X}$, how
one can construct a good black box for the normal closure
$$Y = \left\langle y_1^X,\ldots, y_k^X \right\rangle?$$ One
possibility is to run a random walk on the Cayley graph for $Y$
with respect to the union of the conjugacy classes $S= y_1^X \cup
\cdots\cup y_k^X$ as the generating set for $Y$. If we know that
$Y$ is a simple group then a result of Liebeck and Shalev
\cite[Corollary~1.12]{LiSh} asserts that, for a finite simple
group $G$ and a normal subset $S \subset G$, the diameter on the
Cayley graph ${\mathsf C}(G,S)$ is at most $c\log|G|/\log|S|$. It
might be seen that this result extends to extensions of Lie type
groups by diagonal automorphism and becomes applicable under
conclusions of Fact~\ref{fact:centralisers-of-involutions}
(Shalev, a private communication). It follows from the result by
Diaconis and Saloff-Coste on the mixing time of a random walk on
an edge-transitive graph \cite{DSC} that the mixing time of the
random walk on  ${\mathsf C}(G,S)$ is at most
$c\log^3|G|/\log^2|S|$.

However, we wish to discuss a modification of a product
replacement algorithm whose practical performance as a black box
oracle for normal subgroups seems to be better than a random walk
on ${\mathsf C}(Y,y_1^X \cup \cdots\cup y_k^X)$.

\subsection{Andrews--Curtis graph and the Andrews-Curtis Algorithm}
\label{subsec:AC}

If $G$ is a group (not necessary finite) and  $N \lhd G$, define
the {\em Andrews--Curtis graph\/} $\Delta_k(G,N)$  as the graph
whose  vertices are $k$-tuples of elements in $N$ which generate
$N$ as a normal subgroup:
$$ \Delta_k(G,N) = \left\{\, (h_1,\ldots, h_k) \mid \langle
h_1^G,\ldots, h_k^G\rangle = N \,\right\}.$$ Of course, if the
group $N$ is simple then the vertices of $\Delta_k(G,N)$ are all
$k$-tuples in $N^k \smallsetminus \{(1,\ldots, 1)\}$. Two
vertices are connected by an edge if one of them is obtained from
another by one of the moves
\begin{eqnarray*}
(x_1, \ldots, x_k) & \longrightarrow & (x_1, \ldots,x_ix_j^{\pm 1}, \ldots,  x_k), \; i \ne j \\[2ex]
(x_1, \ldots, x_k) & \longrightarrow & (x_1, \ldots,x_j^{\pm 1}x_i, \ldots,  x_k), \; i \ne j \\[2ex]
(x_1, \ldots, x_k) & \longrightarrow &  (x_1,
\ldots,x_i(x_j^w)^{\pm 1}, \ldots,  x_k), \; i \ne j,
\; w\in G\\[2ex]
(x_1, \ldots, x_k) & \longrightarrow &  (x_1, \ldots,(x_j^w)^{\pm
1}x_i, \ldots,  x_k).
\end{eqnarray*}

Notice that the moves are invertible and thus give rise to a
non-oriented graph.

\begin{conjecture} \textsf{A random walk on the Andrews--Curtis graph
$\Delta_k(G,N)$ provides a `good' black box for $N$. }\end{conjecture}

In practice, a modification of the process, when the last changed
component of the generating tuple (say, $x_ix_j^{\pm 1}$) is
multiplied into the cumulative product $x$, appears to be more
effective:

\begin{itemize} \item \textsc{initialise} $x :=1$. \item \textsc{repeat}
\begin{itemize} \item[$\circ$] Select random $i \ne j$ in $\{\,1,\ldots, k\,\}$.
\item[$\circ$]
\begin{itemize}
\item With equal probabilities, replace $x_i := x_ix_j^{\pm 1}$ or $x_i := x_j^{\pm 1}x_i$,
or \item produce random $w\in G$ and replace
$$x_i := x_i(x_j^w)^{\pm 1} \hbox{ or } x_i := (x_j^w)^{\pm 1}x_i.$$
\end{itemize} \item[$\circ$] Multiply $x_i$ into $x$: $$x := x\cdot x_i.$$ \end{itemize}

\item Use $x$ as the running \textsc{output} of a black box for $N$.
\end{itemize}

Using results on Markov chains, Leedham-Green and O'Brien
\cite{L-GO'B} had shown that  the distribution of values of the
cumulative product $A$ converges exponentially to the uniform
distribution on $N$. However, the issue of explicit estimates is
open and represents a formidable problem.

A discussion of some
related computer experiments can be found in \cite{bl-gnf} and \cite{acbb}.

\subsection{The Andrews--Curtis Problem}

Virtually nothing is known about the properties of the
Andrews--Curtis graph for the free and relatively free groups.
This is one of the few positive results:

\begin{fact}[{\rm A.~G.~Myasnikov \cite{My}}]
For the free solvable group $F_n^{(m)}$ of class $m$ and all\/ $k
\geqslant n$, the Andrews--Curtis graph\/ $\Delta_k(F_n^{(m)},
F_n^{(m)})$ is connected.  \label{fact:myasnikov}
\end{fact}

However, the landscape is dominated by the  Andrews--Curtis
Problem (1965):

\begin{problem}[{\rm Andrews and Curtis \cite{AC}}]
\textsf{Is it true that, for\/ $k \geqslant 2$, the Andrews--Curtis graph
$\Delta_k(F_k, F_k)$ is connected?}
\end{problem}

There is an extensive literature on the subject, see for example,
\cite{AK,BuM,HM}. Some potential counterexamples (originating in
topology) are killed by application of {\em genetic algorithms\/}
\cite{My1}. For example, contrary to the suggestion made by
Akbulut and Kirbi  \cite{AK} in 1985, the pairs  $(x^2y^{-3},
xyxy^{-1}x^{-1}y^{-1})$ and $(x,y)$ of elements in the
2-generator free group $F_2 = \langle x, y\rangle$ belong to the
same connected component of $\Delta_2(F_2, F_2)$.

The work \cite{acbb} suggests a possible line of attack at this
problem based on the study of the connected components of the
Andrews--Curtis graphs
$\Delta_k(G,G)$ for finite groups $G$.

\section{Factor groups}

Assume that we are given a black box group ${X}$ and its normal
black box subgroup ${Y}$. The computations in the factor group
${X}/{Y}$ require testing when two elements $u$ and $v$ in ${X}$
are equal in the factor group ${X}/{Y}$, which is equivalent to
the \textit{membership problem} for ${Y}$:
\begin{quotation}
Check, in polynomial time of $\log |{X}|$, whether $uv^{-1} \in
{Y}$.
\end{quotation}

If ${Y}$ is simple and we have an order oracle for ${X}$, then
the following simple and beautiful algorithm due to Leedham-Green
resolves the membership problem in polynomial time.

\medskip
\noindent \textsc{input:} an element $u \in {X}$.
\begin{itemize}
\item \textsc{for} sufficiently many random
$y_1,\ldots, y_k \in  {Y}$ compute
$$
D := {\rm gcd}(o(uy_1),\ldots, o(uy_k)).
$$
\item \textsc{if} $D = 1$ \textsc{return} $u \in
{Y}$

\textsc{else} \textsc{return} ``probably $u \not\in {Y}$''.
\end{itemize}

This is a one-sided algorithm: if $D =1$ then $u$ definitely
belongs to ${Y}$ for otherwise $D$ is divisible by  the order of
the element $u$ in the factor group ${X}/{Y}$. On the other hand,
orders of sufficiently many random elements of a simple group are
likely to have no non-trivial divisors in common \cite{bps}. See
\cite[Section~4.4]{bs} for a detailed discussion.

\section{Centralisers of involutions}
\label{sec:centralisers}

It is well known that if $u$ and $v$ are involutions in a finite
group, then the group $\langle x,y\rangle$ is a dihedral group;
indeed,
$$
(uv)^u = u^{-1} \cdot uv \cdot u = vu = (uv)^{-1}
$$
and similarly $(uv)^v = (uv)^{-1}$. Hence $u$ and $v$ invert every
element in the cyclic group $\langle uv\rangle$. If the element
$uv$ is of even order then $u$ and $v$ invert the involution
${\rm i}(uv)\in \langle uv \rangle$ and centralise it. If,
however, the element $uv$ has odd order then, by the
Sylow Theorem, the involutions $u$ and $v$ are conjugate by an
element from $\langle uv \rangle$.

This simple observation, due to Richard Brauer, was the starting
point of his programme of classification of finite simple groups
in terms of centralisers of involutions. Remarkably, in the
context of black box groups it can be developed into an efficient
algorithm for constructing black boxes for the centralisers of
involutions.

Let ${X}$ be an arbitrary black box finite group and assume that
$x^E =1$ for all elements $x\in {{X}}$. Write $E = 2^t \cdot r$
with $r$ odd. Let $x$ be a random element in ${X}$. Notice that

\begin{itemize}
\item if $x$ is of odd order, then $x^r=1$ and
$y = x^{(r+1)/2}$ is a square root of $x$:
$$
y^2 =  x^{(r+1)/2} \cdot x^{(r+1)/2} = x^{r+1} = x;
$$

\item if $x$ is of even order then $x^r$ is a $2$-element
and the consecutive squaring of $x^r$ produces the
involution ${\rm i}(x)$ from the cyclic group $\langle x \rangle$.
\end{itemize}

Furthermore, the elements $y = \sqrt{x}$ and ${\rm i}(x)$ can be
found by $O(\log E)$ multiplications.

Let now $i$ be an involution in ${X}$. Construct a random element
$x$ of ${X}$ and consider $z = ii^x$.

\begin{itemize}
\item If $z$ is of odd order and $y = \sqrt{z}$ then
$$
i^y = y^{-1}i y = i yy = iz = i\cdot ii^x = i^x
$$
and $yx^{-1} \in C_X(i)$. We write $ yx^{-1} = \zeta_1(x)$.

\item if $z$ is of even order then ${\rm i}(z)$
lies in the center of the dihedral group $\langle i, i^x\rangle$ and thus ${\rm i}(z) \in C_X(i)$.
We write ${\rm i}(z) = \zeta_0(x)$.

\end{itemize}
Notice that $\zeta_1(x)$ can be computed without knowing the
order $o(x)$ of $x$. One can test whether an element has odd order
by raising it to the odd part $r$ of $E$, and if $o(x)$ is odd
then $x^{(r+1)/2} = x^{(o(x)+1)/2}$.

\medskip

Thus we have a map $\zeta = \zeta_1 \sqcup \zeta_0$ defined by
\begin{eqnarray*}
\zeta: X & \longrightarrow &  C_X(i)\\
x & \mapsto & \left\{ \begin{array}{ll}
\zeta_1(x) = (ii^x)^{(r+1)/2}\cdot x^{-1} & \hbox{ if } o(ii^x) \hbox{ is odd}\\
\zeta_0(x) = {\rm i}(ii^x)  &  \hbox{ if } o(ii^x) \hbox{ is even.}
\end{array}\right.
\end{eqnarray*}

If $c \in G_X(i)$ then
\begin{eqnarray*}
\zeta_0(xc) & =&  {\rm i}(i\cdot i^{xc}) = {\rm i}(i^c\cdot i^{xc}) =
{\rm i}((i\cdot i^x)^c) = {\rm i}(ii^x)^c\\
& = & \zeta_0(x)^c,\\[2ex]
\zeta_1(cx) & = & (ii^{cx})^{(r+1)/2}\cdot x^{-1}c^{-1} =
(ii^{x})^{(r+1)/2}\cdot x^{-1}c^{-1}\\
& = &\zeta_1(x)\cdot c^{-1}.
\end{eqnarray*}

Hence if the elements $x \in {X}$ are uniformly distributed and
independent in $X$ then
\begin{itemize}
\item  the distribution of elements $\zeta_1(x)$ in $C_X(i)$
is invariant under  right multiplication by elements $c \in C_X(i)$,
that is, if $A \subset  C_X(i)$ and $c \in C_X(i)$ is an arbitrary element then
the probabilities $P(\zeta_1(x) \in A)$ and $P(\zeta_1(x) \in Ac)$ coincide.

\item The distribution of involutions $\zeta_0(x)$ is
invariant under the
action of $C_X(i)$ on itself by conjugation, that is,
$$
P(\zeta_1(x) \in A) =P(\zeta_1(x) \in A^c).
$$
\end{itemize}

\subsection{Running the odd type oracle}

Therefore we came to the following simple but important result.

\begin{theorem}
If the elements\/ $x \in {X}$ are uniformly distributed and
independent in\/ ${X}$ then the elements\/ $\zeta_1(x)$ are
uniformly distributed and independent in\/ $C_X(i)$.
\end{theorem}

This gives us a good black box for $C_X(i)$; we shall call it the
\textit{ black box of odd type} or \textit{odd type oracle}.

It might happen, however, that the share of elements $x\in {X}$
for which the function $\zeta_1(x)$ is defined is too small to
use $\zeta_1$ as an efficient way to generate elements in
$C_X(i)$. For example, in the group $Y=PSL_2({\mathbb{F}}_{q})$,
when $q= p^k$ is a big power of an odd prime integer $p$, almost
every element is semisimple and thus belongs to a cyclic group of
order $(q\pm 1)/2$; one of these two subgroups has even  order
and at least $1/2$ of its elements are also of odd order. This
shows that between $1/4$ and $1/2$ elements in $Y$ are of even
order. All involutions in $Y$ are conjugate. It is easy to see
that the product of two random involutions in $Y$ has even order
with probability between $1/4$ and $1/2$. When we work in the
direct product $X = Y\times\cdots \times Y$ ($k$ times), we have
to make these computations componentwise, which leads to the
unfortunate conclusion that, in the worst case scenario, the
probability of the product of two conjugate involutions to be of
odd order could be close to $1/2^k$.

If $X$ is a simple group of Lie type of odd characteristic, then
we have the following crude estimate.
\begin{theorem}
Let\/ $G$ be a simple group of Lie type of odd characteristic and
Lie rank $n$. If\/ $i$ is an arbitrary involution in $X$ then the
product\/ $i\cdot i^g$ has odd order with probability $c/n^d$ for
some constants $c$ and\/ $d$.
\end{theorem}
This can be deduced from the Galois cohomology  of {\em reflexive
tori} in simple algebraic groups $G$, that is, tori $T$ such that
$t^i = t^{-1}$ for some involution $i\in N_G(T)$ and all $t \in
T$. This theory is developed in \cite{reflexive} by analogy with
the classical theory of tori in semisimple algebraic groups over
finite fields \cite{ste1}.

The situation is better in simple groups of Lie type over big
finite fields of characteristic $2$, where almost all elements are
semisimple and thus have odd order. In this context, a product of
an involution and its conjugate almost always has  odd order, and
 the odd type oracle works with the maximal possible
efficiency.

 Of course, when we deal with
the verification of a possible isomorphism between a black box
group ${X}$ and the known target group ${G}$, we can try to
locate in ${X}$ an involution $i$ which, in the case of
isomorphism ${X} \simeq {G}$, should behave like an involution
$j$ from a conjugacy class in ${G}$ where the share of elements
$g \in {G}$ with $o(jj^g)$ odd is sufficiently big; then we can
run the odd type black box in the hope to eventually get a
contradiction with the isomorphism ${X} \simeq {G}$. Examples of
this type of computation can be found in \cite{altseimer-borovik,borovik}
and are used in our proof of Theorem~\ref{th:Op(X)}. The {\em
classical involutions} in simple groups $X$ of Lie type of
characteristic $p>3$ are particularly useful. Recall that a
classical involution $t$ has the property that $C_X(t)$ contains
a subnormal subgroup $L \simeq SL_2({\mathbb{F}}_{p^k})$ such
that $t \in Z(L)$. For example, in $X = SL_n({\mathbb{F}}_{p^k})$
a classical involution is an involution with exactly $2$ or $n-2$
eigenvalues $-1$. Classical involutions played the very prominent
role in the classification of finite simple groups \cite{asch},
and their reappearance in the theory of black box groups is not
really surprising.

\subsection{The oracle of even type}

If the function $\zeta_1(x)$ is rarely defined then  the function
$\zeta_0(x)$ is defined for almost all $x\in {X}$ and produces a
normal set of involutions with probability distribution invariant
under conjugation by $C_X(i)$. Hence $C^\circ_X(i) = \langle
\zeta_0(x) \mid x \in X\rangle$ is a normal subgroup in $C_X(i)$.
It seems reasonable to take  the cumulative product $z:=
z\zeta_0(x)$ of consecutive values of $\zeta_0(x)$ for the output
of a black box for $C^\circ_X(i)$. We shall call it the
\textit{oracle of even type}.

The values of the function $\zeta_0$ belong to the union of
conjugacy classes
$$ S = s_1^{C^\circ_X(i)} \cup \cdots \cup s_k^{C^\circ_X(i)}$$
with the probability distribution  invariant under the
conjugation by elements from $C_X(i)$. In the case when $X$ is a
reductive group of characteristic $p>2$,  the analysis similar to
that of Section~\ref{sec:normal} shows that the mixing time
$t_{\rm mix}$ of the corresponding random walk on
$\mathsf{C}(C^\circ_X(i),S)$ is bounded by
$$
t_{\rm mix} < c \max_i \frac{\log^3 |C^\circ_X(i)|}{\log^2
\left|s_i^{C^\circ_X(i)}\right|}.
$$

\subsection{Centraliser of a transposition in a symmetric group}
\label{sec:transposition}
 In one case, namely, when $X = {\rm
Sym}_{n+2}$ and $i = (12)$, the detailed analysis of the even
type oracle is already contained in Diaconis and Shahshahani
\cite{ds}. Indeed, it is easy to see that $C_X(i) \simeq
{\mathbb{Z}}_2 \times {\rm Sym}_{n}$. Random products $ii^x$ are
either elements of order $1$ or $3$ (which happens with
probability $O(1/n)$ when $1^x \in \{1,2\}$ or $2^x \in \{1,2\}$)
or one of the $n \choose 2$ involutions $\zeta_0(x)=(12)(st)$.
Hence the involutions $\zeta_0(x)$ generate a subgroup $C^\circ$
of index $2$ in $C_X(i)$. It is easy to see that $C^\circ$ is
isomorphic to the symmetric group ${\rm Sym}_n$. We see that the
even type oracle $\zeta_0$ works with sufficient speed. As for
the distribution of the images of the cumulative product
$z:=z\cdot\zeta_0(x)$ in the factor group $C^\circ \simeq {\rm
Sym}_n$, we can use an estimate from \cite{ds} regarding
generation of ${\rm Sym}_n$ by random transpositions:
\begin{fact} {\rm (Diaconis and Shahshahani \cite{ds})}
If\/ $k = \frac{1}{2}\cdot n\log n  + cn$, $c >0$, and $P^{*k}$
is the distribution of the random product of $k$ transpositions
then
$$
\| P^{*k} - U\| \leqslant ae^{-2c}
$$
for an universal constant\/ $a$. Here $U$ is the uniform
distribution on ${\rm Sym}_n$ and the norm $\|\,\|$ is defined as
$$
\| P^{*k} -U \| = \frac{1}{2} \sum_{g \in {\rm Sym}_n}
\left|P^{*k}(g)-U(g)\right|. $$
\end{fact}
This means that we have to skip first $\frac{1}{2}\cdot n\log n$
values of the cumulative product $z$, and after that we can
expect that the cumulative product quickly converges to the
uniform distribution on $C$. Further results in \cite{ds} show
that the threshold estimate $\frac{1}{2}\cdot n\log n$ cannot be
improved.

\subsection{The mixed type oracle}

Some computer experiments suggest that the cumulative product
$$
z:= \left\{\begin{array}{rl} z\cdot\zeta_0(x) & \hbox{ if }
o(ii^x) \hbox{ is even}\cr z\cdot \zeta_1(x) & \hbox{ if } o(ii^x)
\hbox{ is odd} \end{array} \right.
$$
has a satisfactory performance as a black box for $C_GX(i)$. We
shall call it the \textit{mixed type oracle}.

\section{Improving black boxes by  cumulative
product} \label{sec:Fourier}

It is obvious that the probability distributions of the outputs
of black boxes as they  appear in practical computations might
considerably deviate from the uniform distribution. Independence
of consecutive values of the output is also questionable (it is
the case, for example, with the product replacement algorithm).
For that reason methods of improving the statistical properties
of the output are highly desirable---of course, if they are
computationally efficient.

Assume that the output $\{\, x_1, x_2, \ldots \,\}$ of our black
box ${X}$ is independent and has the probability distribution $P=
\{p(x)\}$. Then the cumulative product $x:= x \cdot x_i$ is the
random walk on the group $X$ generated by the probabilistic
distribution $P$, that is, a random walk in which we move from
the element $x$ to $xy$ with probability $p(y)$. The probability
distribution after $k$ steps of the random walk is the $k$-th
convolution $P^{*(k-1)}$ of $P$. Here the \textit{convolutions}
are defined as
$$(P*Q)(x) = \sum_{y\in X} P(xy^{-1})Q(y), \quad  P^{*k} = P *  P^{*(k-1)}.
$$

\subsection{Centralisers of involutions}

When the distribution $P$ is invariant  under conjugation by
elements of $X$ (as it happens in generation of centralisers of
involutions), it appears that the methods based on the
non-commutative Fourier transform \cite{diaconis} are useful for
the estimating the rate of convergence of $P^{*(k)}$ to the
uniform distribution on $X$. In particular, under these
conditions the Upper Bound Lemma by Diaconis and Shahshahani
\cite{ds} takes the following form:
$$
\|P^{*(k)} - U\|^2 \leqslant \frac{1}{4}\sum_{\chi \ne 1}\left(
\chi(1)^2\left|\sum_x P\left(x^G\right)
\frac{\chi(x)}{\chi(1)}\right|^{2k}\right),
$$
where  the first sum is taken over the non-trivial irreducible
characters of $G$ and the second over representatives of
conjugacy classes of $X$. This is how  Diaconis and Shahshahani
obtained their bounds for the generating of ${\rm Sym}_n$ by
random transvection \cite{ds} used in
Section~\ref{sec:transposition} in the analysis of generation of
the centraliser of a transposition.

It would be interesting to get some numerical data related to
generation of centralisers of involutions. The first step of this
problem seems to be relatively easy.
\begin{question}
\textsf{ For  an involution $i$ in a finite simple group $X$ of Lie type of
odd characteristic, compute the probability distribution $P(z)$
of the values of the function $\zeta_0: X \longrightarrow
C^\circ_X(i)$.}
\end{question}
However, obtaining explicit estimates in the Upper Bound Lemma
requires a detailed knowledge of characters of the group
$C^\circ_X(i)$. It would be interesting to complete the analysis
at least in some number of `small' cases.

\subsection{The Andrews--Curtis algorithm} Notice that if we slightly
modify the Andrews--Curtis algorithm and conjugate the both
elements involved in multiplication,
$$
x_i := x_i^u(x_j^v)^{\pm 1} \hbox{ or } x_i := (x_j^v)^{\pm
1}x_i^u, \quad i \ne j
$$
where $u$ and $v$ are random elements of $X$, then the
distribution of the output of the new algorithm becomes invariant
under the action of $X$ by conjugation. Assume that $Y$ is
simple. If we make the assumption that the consecutive values of
the output are independent (which is sufficiently close to the
truth when the size $k$ of  generating $k$-tuples is sufficiently
big), the combination of results of Liebeck and Shalev (see
Section~\ref{sec:normal}) with the results of Diaconis and
Saloff-Coste \cite{DSC} gives a cubic (in $\log |G|$) estimate
for the mixing time of the cumulative product. The experimental
data shows  a much better performance of the Andrews--Curtis
algorithm.  It would be very interesting to carry out the
rigorous analysis of the behaviour of the cumulative product in
the Andrews--Curtis algorithm.

\subsection{The product replacement algorithm} As it was demonstrated
by Leedham-Green et al.\ \cite{bl-gnf,L-GO'B}, the practical
performance of the product replacement algorithm can be improved
if we multiply its consecutive outputs $z_1, z_2,\ldots$ into the
cumulative product
 $z := z \cdot z_i$, and take the consecutive values of  $z$ for
the output of a new black box. We noticed that the cumulative
product gave a similar improvement in the black box algorithm for
centralisers of involutions. It would be interesting to formulate
and prove results which confirm or at least provide a heuristic
justification of these empirical observations.

\bigskip

\centerline{\textsc{Acknowledgements}}

\medskip

The author thanks Bill Kantor, Charles Leedham-Green, Martin
Liebeck, Eamonn O'Brien, Igor Pak,  Jan Saxl, and Aner Shalev for
helpful discussions.

\vfill

\noindent
 Department of Mathematics,  UMIST, PO Box 88,
 Manchester M60 1QD, UK
 \\
 \textsc{E-mail:} \texttt{borovik@umist.ac.uk}\\
\textsc{Web:} \texttt{http://www.ma.umist.ac.uk/avb/}
 \end{document}